\documentclass[12pt,leqno]{article}

\usepackage{amscd}
\usepackage{amssymb}
\usepackage[dvips]{graphics}
\usepackage{epic,eepic,epsfig}
\usepackage{amsmath}

\usepackage{amsfonts}
\usepackage[all]{xy}

\def\C{\mathbb{C}}
\def\H{\mathbb{H}}

\def\P{\mathbb{P}}

\def\Q{\mathbb{Q}}
\def\R{\mathbb{R}}
\def\Z{\mathbb{Z}}

\def\A{\mathcal{A}}

\def\D{\mathcal{D}}

\def\F{\mathcal{F}}

\def\L{\mathcal{L}}

\def\O{\mathcal{O}}

\def\codim{\text{\rm codim }}
\def\dim{\text{\rm dim }}

\def\ker{\text{\rm Ker }}
\def\rank{\text{\rm rank }}

\newtheorem{thm}{Theorem}[section]
\newtheorem{cor}{Corollary}[section]

\newtheorem{example}{Example}[section]

\newtheorem{remark}{Remark}[section]

\begin{document}

\pagestyle{myheadings}
\markboth{}{}

\title{\bf  Local topology of reducible divisors }

\author{\bf Alexandru Dimca and Anatoly Libgober}

\date{}

\maketitle

\begin{abstract} We show that the universal abelian 
cover of the complement to a germ of a reducible 
divisor on a complex space $Y$ with isolated singularity 
is $({\rm dim} Y-2)$-connected provided that the divisor 
has normal crossings outside of the singularity of $Y$.
We apply this result to obtain a vanishing property for the 
cohomology of local systems of rank one and also 
study vanishing in the case of local systems of 
higher rank.

\end{abstract}

\section{Introduction}

\bigskip \

The topology of holomorphic functions near an isolated
singular point is a classical subject (cf. \cite{Milnor},
\cite{D1}). Among the main
results are the existence of Milnor fibration 
and the connectivity of the Milnor fiber yielding
a very simple picture for the latter: it has the homotopy
type of a wedge of spheres. Starting with the case 
of a germ of holomorphic function on 
${\C}^N$ considered by Milnor
(\cite{Milnor}), these results were eventually extended to
the germs of holomorphic functions on analytic spaces
(cf. \cite{H1},  \cite{Le2}  ).\\ 

In \cite{INNC}, it was shown 
that if the divisor of a holomorphic
function on ${\C}^N$ is {\it reducible} then these
connectivity results can be refined. This refinement is 
based on an 
observation that the Milnor fiber
is homotopy equivalent to the infinte cyclic cover
of the total space of the Milnor fibration. In the case
when the divisor of a holomorphic function is reducible
one has the associated abelian cover which has interesting
connectivity properties generalizing the connectivity properties in 
the cyclic case. The present paper studies the case of reducible 
divisors on {\it arbitrary} analytic spaces with isolated singularities.\\

More precisely the situation we consider is the following.
Let $(Y,0) \subset (\C^N,0)$ be an $(n+1)$-dimensional complex analytic space germ with an isolated singularity at the origin. Let $(D_j,0) \subset (Y,0)$ for ${j=1,...,r}$ be $r$ irreducible Cartier divisors on $(Y,0)$.
 We set $X=\cup_{i=1,r} D_i$, $M=Y \setminus X$ and regard $M$ as the complement of the hypersurface arrangement $\D=(D_j)_{j=1,r}$.\\

In this paper we investigate the topology of this complement $M$. In Section 2 we generalize a case of the L\^e-Saito result in \cite{LeSaito} 
asserting that if $(Y,0)$ is a smooth germ and $X$ is an isolated non-normal crossing divisor (see the definition below), then the fundamental group $\pi _1(M)$ is abelian. Our proof is based on an idea used in \cite{Nori} 
in global case and is much shorter than the proof in  \cite{LeSaito}.\\

In Section 3 we consider the case when the  hypersurface arrangement  $\D$
is an arrangement based on a hyperplane arrangement $\A$ in the sense of Damon
\cite{Da}. We show that the (co)homology of $M$ is determined up-to degree $(n-1)$ by the hyperplane arrangement $\A$. The key fact here is the functoriality
of the Gysin sequence and the splitting of the  
Gysin sequence associated to a triple $(\A,\A',\A'')$ of hyperplane arrangements into short exact sequences.
Combining the results above and following the approach in \cite{INNC}, we show that the universal abelian covering $\tilde M$ of $M$ is homotopically a bouquet of spheres of dimension $n$ which is the refinement of 
\cite{Milnor} and \cite{H1} we mentioned earlier.\\

In the last two sections we prove vanishing results for the (co)homology of
the complement $M$ with coefficients in a local system $\L$ on $M$.
The case when the $\rank$ of $\L$ is equal one is treated in Section 4 and in this 
context we give a description for the dimension of the non zero homology groups $H_*(M,\L)$.
The general case  when $\rank \L \geq 1$ is treated 
in Section 5 where we allow a
more general setting for the ambient space $(Y,0)$ and for the divisor $(X,0)$.
The vanishing result in this case follows the general philosophy in \cite{EV},
but the use of perverse sheaves as in \cite{CDO} is unavoidable. Note that in our case the space $M$ may be singular so one cannot use the technique of integrable connections to get vanishing results. 
A new point in our proof is the need to use the interplay between constructible
complexes of sheaves on real and complex spaces. Indeed, real spaces occur in the picture in the form of links of singularities.

\section{Fundamental group of the complements to INNC.}

Let $(Y,0) \subset (\C^N,0)$ be as above an $(n+1)$-dimensional complex analytic space germ with an isolated singularity at the origin. Let $(D_j,0) \subset (Y,0)$ for ${j=1,...,r}$ be $r$ irreducible Cartier divisors on $(Y,0)$, i.e. each $D_j$ is given (with its reduced structure) as the zero set of a holomorphic function germ $f_j:(Y,0) \to (\C,0)$.
When the local ring ${\cal O}_{Y,0}$ is factorial, then any hypersurface germ in $(Y,0)$ is Cartier. This is the case for instance when  $(Y,0)$ is an isolated complete intersection singularity (ICIS for short) with $\dim Y \geq 4$, see
\cite{Gr}. See also Example \ref{ex1}.
 Here and in the sequel we identify germs with their (good) representatives.\\ 
 We assume in this section that the following condition holds.

%\bigskip (*) The link $L_Y$ of $Y$ at  
%$0$, i.e. the intersection of $Y$ with 
%the boundary $\partial B_{\epsilon}$ 
%of a small ball $B_{\epsilon}$ about $0$, is simply connected.

\bigskip 

(${\bf C1}$) The divisor $X=\cup_{i=1}^{i=r}D_i$ 
has only normal crossing singularities on $Y$
except possibly at the origin. We say in this case that $X$ is an isolated non normal crossing divisor (for short INNC) on $(Y,0)$.

\bigskip 

In particular each germ $(D_j,0)$ has an 
isolated singularity at the origin as well. 
Since the $(r+1)$-tuple $(Y,D_1,...,D_r)$ has a conical structure 
(cf. \cite{GM1})
we have an isomorphism:

\begin{equation}
\pi_1(L_Y \setminus  L_X)=\pi_1(M \cap \partial B_{\epsilon}) \rightarrow 
\pi_1 (M)
\end{equation}
where $L_Y$ (resp. $L_X$) denotes the link of $Y$ (resp. of $X$), i.e.
 the intersection of $Y$ with 
the boundary $\partial B_{\epsilon}$ 
of a small ball $B_{\epsilon}$ about $0$. In particular we get a morphism
$$\pi_1(M)=\pi_1(L_Y \setminus  L_X) \to \pi_1(L_Y)$$
induced by the inclusion $L_Y \setminus  L_X \to L_Y$.

\begin{thm}\label{abelianfg}
 For $n \ge 2$, the kernel of the surjection
 $\pi_1(M) \rightarrow \pi_1(L_Y)$ is contained in the center
of the group $\pi_1(M)$. In particular if $L_Y$ is simply connected
then the fundamental group $\pi_1 (M)$ is abelian.
\end{thm}

\begin{proof} 

First notice that if ${\rm dim} Y >3$ and if $H$ 
is a generic linear subspace 
passing through $0$ such that  
the codimension of $H$ in ${\C}^N$
is ${\rm dim} Y-3$, 
then, by Lefschetz hyperplane section theorem (cf. \cite{GM1}, p. 26 and p. 155),
we have an isomorphism:
\begin{equation}
\pi_1 (M \cap H) \rightarrow \pi_1 (M)
\end{equation}
Hence it is enough to consider the case ${\rm dim} Y=3$ only
(though the arguments below work for any dimension).

Next notice that $\kappa={\rm Ker}(\pi_1(L_Y \setminus  L_X) \rightarrow 
\pi_1(L_Y))$ is 
the normal closure of the set of elements in the fundamental group 
$\pi=\pi_1(L_Y \setminus L_X)$ 
represented by the loops $\delta_i$ each of which is the  
boundary of a fiber  of a small tubular neighbourhoods 
of $D_i \cap L_Y$ in $L_Y$. 
Indeed a loop representing an element
$\gamma$ in the kernel $\kappa$
is the image of the boundary of a 2-disk under
a map $\phi: D^2 \rightarrow L_Y$ which is isotopic to an
embedding (since ${\rm dim}L_Y=5$) and which we may 
assume to be transversal to all $D_i \cap L_Y$. Now $\delta_i$ 
are the $\phi$-images of loops in $D^2$ each of which is
composed of a path $\alpha_i$ going from the point in $D^2$ 
corresponding to the 
base point $p \in L_Y$ to the vicinity of points corresponding 
to $\phi (D^2) \cap D_i$, a small loop about $\phi(D^2) \cap D_i$ 
and back along $\alpha_i^{-1}$. 
So it is enough to show that all $\delta_i$ belong
to the center of $\pi_1(L_Y \setminus L_X)$. 

Let $T(D_i)$ be a tubular neighbourhood of $D_i \cap L_Y$ in $L_Y$
as above.
We claim that for any $i$ ($i=1,..,r$) there is a surjection:
\begin{equation}\label{surjection}
\pi_1(T(D_i) \setminus X) \rightarrow \pi(L_Y \setminus L_X)
\end{equation}
Notice that assuming (\ref{surjection}) we can conclude the proof as follows. 
Since $D_i$'s have normal intersections in $L_Y$, the space
$T(D_i) \setminus X   $ 
is homotopy equivalent to the total space of a locally trivial circle 
fibration over $D_i \setminus \cup_{j \ne i}D_j$.
The fiber $\delta_i^{\prime}$ of this fibration, which is a loop based at a point
$p^{\prime}$, is in the center of 
$\pi_1(T(D_i) \setminus X   ,p^{\prime})$. Therefore 
(\ref{surjection}) yields that image of $\delta^{\prime}$ commutes
with any element in $\pi_1(L_Y \setminus L_X  ,p^{\prime})$ and hence
with any element in $\pi_1(L_Y  \setminus L_X,p)$.

To show the surjectivity (\ref{surjection}), let us consider a sufficiently generic 
holomorphic function $g$ on $Y$ 
with the set of zeros sufficiently close to the zero set of $f_i$
so that $L_Y \cap \{g=0\} \subset T(D_i)$.
We have the decomposition 
\begin{equation}
\begin{matrix} \pi_1 (L_Y \cap \{g=0\} \setminus L_X ) &  & \cr
              \downarrow & \searrow & \cr    
            \pi_1 (T(D_i)\setminus L_X   ) & \rightarrow  & \pi_1 (L_Y\setminus L_X) \cr
\end{matrix}
\end{equation}
corresponding to the factorization 
of the embeddings. This yields that the horizontal map is surjective
provided the map:
\begin{equation}\label{gsurjects}  
\pi_1(L_Y \cap (g=0)\setminus L_X) \rightarrow \pi_1(L_Y\setminus L_X)
\end{equation}  
is surjective. But this follows from \cite{LeHamm87}.
%This can be obtained as an application of Morse theory.
%We can replace $L_Y$ by the intersection of $Y$ with the 
%boundary of polydisk $\vert u_i \vert \le \epsilon, i=1,...,N$.
%The boundary of this polydisk consists of $N$ polydisks
%of dimension $N-1$ given by $\vert u_i \vert =\epsilon, 
%\vert u_j \vert \le \epsilon$. Since the codimension of the space
%of functions $g$ which are tangent to a codimension 2 subspace 
%in ${\C}^N$ given by $u_i=u_j=0$ has codimension at least 2,
%we can assume that there are no critical points of $\vert g \vert $ 
%on the edges of the polydisk. For a sufficiently generic $g$ 
%the critical points of $\vert g \vert $ restricted on 
%$(\vert u_i \vert=\epsilon) \cap Y$ correspond to the 
%singularities of $g$ on $u_i=\epsilon$ and have index 2
%which yields the surjectivity of (\ref{gsurjects}).

\end{proof}

\bigskip

Note that this result in the case when $Y={\C}^{n+1}$ 
is a consequence of a theorem of L\^e Dung Trang and K.Saito
(cf. \cite{LeSaito}).

\begin{example} \label{ex1} \rm

(i) If $(Y,0)$ is an ICIS with $\dim (Y,0) \geq 3$, then it follows from \cite{H1} that the link $L_Y$ is simply-connected.

(ii) If $V \subset \P^m$ is a locally complete intersection such that
$n=\dim V> \codim V$, then the morphism $\pi_2(V) \to \pi_2( \P^m)=\Z$ induced by the inclusion $V \to \P^m$ is an epimorphism by the generalized Barth Theorem, see \cite{GM1}, p. 27. It follows that the associated affine cone $(Y,0)=(CV,0)$ has a simply-connected link. 

If we assume that $V$ is smooth  and that $n=\dim V> \codim V+1$, then one has $C\ell(\O _{Y,0})=0$,
i.e. any divisor on this germ $(Y,0)$ is Cartier. This follows from the exact sequence in \cite{Ha}, Exercise II.6.3 comparing the divisor class groups in the local and the global settings, the usual isomorphism 
$C\ell(V)=H^1(V,{\cal O}^*_V)$, see \cite{Ha}, II.6.12.1 and II.6.16 and the GAGA results allowing to use the exponential sequence, see \cite{Ha}, Appendix B, to relate topology to $H^1(V,{\cal O}^*_V)$.

If $E=\cup _{j=1,r}E_j$ is a normal crossing divisor on the smooth variety $V$, then the associated cone $X=\cup _{j=1,r}CE_j$
is an INNC divisor on $(Y,0)$.\\
The epimorphism $\pi_1(M) \to \pi_1(V \setminus E)$ can be used to show that this last fundamental group is abelian.

\end{example}

\section{Homology of the complements to reducible divisors}

Assume in this section that the germ $(Y,0)$ is an ICIS with $dimY=n+1$ and let $\A=\{H_i\}_{i=1,...,r}$ be a central hyperplane arrangement in $\C^m$. Suppose given an analytic  map germ
$f:(Y,0) \to (\C^m,0)$ such that the following condition holds.

\bigskip

({\bf C2}): For any edge $L \in L(\A)$ with $\codim L =c$, the (scheme-theoretic)
pull-back $D_L=f^{-1}(L)$ is an ICIS in $(Y,0)$ of codimension exactly $c$ for $c \leq n$ and  $D_L =\{0\}$ for $c \geq n+1$.

\bigskip

This condition $({\bf C2})$ is equivalent to asking that $f:(Y,0) \to (\C^m,0)$ is transverse to $\A$ off $0$ in the sense of Damon, see \cite{Da}, Definition 1.2.
In his language, $X= \cup _{i=1,r}D_i$ is a nonlinear arrangement of hypersurfaces based on the central arrangement $\A$, with $D_i=f^{-1}(H_i)$. Note that
for $n \geq 2$ all the germs$D_j$ are irreducible by \cite{H1}, but on the other hand the condition $({\bf C1})$ may well fail in this setting.

\bigskip

Consider the complements $M=Y \setminus X$ and
$N=\C^m  \setminus \cup  _{i=1,r}H_i$, and note that there is an induced mapping
$f:M \to N$. Our result is the following

\begin{thm}\label{homologicalthm}

With this notation, 
$$f_*:H_j(M) \to H_j(N)$$
is an isomorphism for $j<n$ and an epimorphism for $j=n$. Similarly
$$f^*:H^j(N) \to H^j(M)$$
is an isomorphism for $j<n$ and a monomorphism for $j=n$. In particular, the algebra $H^*(M)$ is spanned by $H^1(M)$ up-to degree $(n-1)$.
\end{thm}

\begin{proof}

For $n=1$ everything is clear, so we can assume in the sequel $n>1$.
 
The proof is by induction on $r$.
 For $r=1$
the result follows since $M$ can be identified to the total space of the Milnor fibration, whose Milnor fiber is a bouquet of $n$-dimensional spheres by work of Hamm, see \cite{H1}.

Assume now that $r>1$ and apply the deletion and restriction trick, see more on this in \cite{OT}. Namely,
let $\A'= \{H_i\}_{i=2,...,r}$ and $\A''= \{H_1 \cap H_i\}_{i=2,...,r}$.
Then both  $\A'$ and  $\A''$ are central arrangements with at most $(r-1)$
hyperplanes.

Since $L(\A') \subset L(\A)$, it is clear that $f:(Y,0) \to (\C^m,0)$
satisfies the condition $(C)$ with respect to $\A'$. Moreover,
 $L(\A'') \subset L(\A)$ (with a 1-shift in codimensions) and
$f:(D_1,0) \to (\C^{m-1},0)$
satisfies the condition $(C)$ with respect to $\A''$.

Let $M',M'',N',N''$ be the corresponding complements. Since $M''$ (resp. $N'')$
is a smooth hypersurface in $M'$ (resp. $N'$) we have the following ladder of Gysin sequences.

$$
\xymatrix{\cdots \ar[r] & H_{j-1}(M'') \ar[d]_{f_*} \ar[r] & H_j(M) \ar[d]_{f_*} \ar[r] &
                            H_j(M') \ar[d]_{f_*}  \ar[r] &  H_{j-2}(M'') \ar[d]_{f_*} \ar[r]    &     \cdots \\
 \cdots \ar[r] &  H_{j-1}(N'')  \ar[r] &  H_j(N)   \ar[r] &  H_j(N')     \ar[r]
& H_{j-2}(N'')   \ar[r]  &  \cdots 
}$$

It is a standard fact in hyperplane arrangement theory that the bottom-right morphism $ H_j(N')  \to  H_{j-2}(N'')$ is zero. Usually this result is stated for cohomology, see \cite{OT}, p.191, but since both homology and cohomology of  hyperplane arrangement complements are torsion free, the vanishing holds for homology as well.

An easy diagram chasing, using the induction hypothesis, shows that for $j \leq n$ the top-right morphism  $ H_j(M')  \to  H_{j-2}(M'')$ is zero as well. This implies the claim, again by an easy diagram chasing and using the induction hypothesis.

The proof for the cohomological result is completely dual.

\end{proof}

Using Theorems \ref{abelianfg} and \ref{homologicalthm} and Example \ref{ex1},(i) we get the following result.

\begin{cor} \label{pi=H}

Assume that $n \geq 2$ and that the divisor $X$ satisfies the condition $(C1)$.Then
$$\pi_1(M)=H_1(M)=H_1(N)$$
is a free abelian group of rank $r=|\A|$.
\end{cor}

\begin{remark} \rm

Note that $M$ is a Stein manifold, and hence $H_j(M)=0$ for $j>n+1$ by \cite{H2}.
It follows that there are only two Betti numbers $b_j(M)$ to compute, 
namely for $j=n,n+1$ since for $j<n$ one has $b_j(M)=b_j(N)$ known by the results in \cite{OT}. 
Moreover by the additivity of Euler characteristics, see \cite{F}, it follows that $\chi(M)= \chi(Y) -\chi(X)= 0$. This gives a relation between the two top unknown Betti
numbers of $M$.

\end{remark}

Similarly to 
\cite{INNC}, the above results yield the following.

\begin{thm} \label{univcover} Let $(Y,0)$ be a  
an isolated complete intersection 
singularity of dimension 
$n+1 \geq 3$.
Let $X=\cup_{i=1,r} D_i$ be a union of Cartier divisors of $Y$ which have normal crossing 
outside of the origin. Then the universal abelian cover ${\tilde M}$ of $M=Y \setminus X$ 
has the homotopy type of a bouquet of spheres of dimension $n$.
\end{thm} 

\begin{proof} The proof is similar to the proof of Thm. 2.2 in \cite{INNC}.
Firstly, let us consider the exact homotopy sequence 
corresponding to the map $f: M \rightarrow {{\C}^*}^r$, obtained by using the equations $f_i=0$ for the divisors $D_i$.
The isomorphism $\pi_1(M) \rightarrow \pi_1({{\C}^*}^r)$ which follows from
 Theorems \ref{abelianfg} and \ref{homologicalthm}
and the known fact $\pi_j({{\C}^*}^r)=0$ for $j>1$ yield that 
$\pi_2({{\C}^*}^r,M)=0$ and $\pi_j({{\C}^*}^r,M)=\pi_{j-1}(M)$ for $j>2$ (here we assume that 
$f$ is replaced by an embedding, which is of course possible up-to homotopy type). Moreover we can show exactly as in \cite{INNC} that the action of $\pi_1(M)$ on  $\pi_j({{\C}^*}^r,M)$ is trivial.
Hence we can apply the relative Hurewich theorem to the 
pair $({{\C}^*}^r,M)$ and note that we have a vanishing 
of the relative homology of this pair as a 
consequence of the previous theorem. Looking now at $f: M \rightarrow {{\C}^*}^r$ as a homotopy fibration with fiber $\tilde M$, we get the vanishing
of the homotopy groups of the universal abelian cover  ${\tilde M}$ of 
$M$ up to dimension $n-1$. On the other hand, the 
existence of the
Milnor fibration of $g=f_1...f_r:(Y,0) \to (\C,0)  $ (theorem of Hamm in \cite{H1}) yields that 
$M$ admits a cyclic cover   which has the homotopy type 
of a CW complex of dimension $n$ (i.e. the Milnor fiber $F$ of 
the hypersurface $X$ in $Y$ ). Hence the universal abelian 
cover  ${\tilde M}$, which is the universal abelian cover of this Milnor fiber $F$ has the 
homotopy type of an $n$-complex. Therefore the universal abelian cover  ${\tilde M}$
is homotopy equivalent to the wedge of spheres $S^{n}$.

\end{proof}

\section{Homology of local systems (rank one case)}

Let $(Y,0)$ be a germ of an isolated complete intersection singularity 
and let $X=\bigcup_{i=1,r} D_i$ be a divisor which has normal crossings
outside of the origin, i.e. we place ourselves
again in the setting of Theorem \ref{univcover}. The above notation is still used here.\\
Let $\rho: \pi_1(M) \rightarrow {\C}^*$ 
be a character of the fundamental 
group or equivalently a local system $\L$ of rank one on $M$. The space
$M$, being a Stein space of dimension $(n+1)$, has the homotopy type of an $(n+1)$ complex 
and hence $ H_j(M,\L) =  H_j(M,\rho)=0$ for $ j >n+1$.
The main result of this section is the following:

\begin{thm}\label{locsys}
 Let $\rho: \pi_1(M) \rightarrow {\C}^* $ be a non trivial character and let
$\L$ be the associated rank one local system on $M$. 
Then:

1. $ H_j(M,\L) =0$ for  $j \ne n,n+1 $;

2. ${\dim} H_n(M,\L)$ is the largest integer $k$
 such that $\rho$ belongs to the zero set $V_k$  of 
the $k$-th Fitting ideal of the ${\C}[\pi_1(M)]$-
module $\pi_n(M) \otimes_{\Z} {\C}$.

3.The largest integer $k$ such that the trivial character of $\pi_1(M)$
belongs to $V_k$ is equal to  
$${\dim \ker} (\Lambda^{n+1}H^1(M) \rightarrow 
H^{n+1}(M))+{\dim} H_n(M)-{r \choose n}$$.

\end{thm}

\begin{proof}

Recall the spectral sequence for the cohomology of 
local systems (cf. \cite{CE}) Thm. 8.4).
Let $C_*^{\rho}({\tilde M})$ be the chain complex on which 
$H_1(M,{\bf Z})$ acts from the right via $g(x)=\rho(g)xg 
\ (g \in H_1(M,{\Z}), x \in C_*(\tilde M,{\C}))$ 
where $x \rightarrow xg$ is the action via deck transformations.
Let $H_q^{\rho}(\tilde M)$ be the homology of this complex.
We have a spectral sequence:
$$E^2_{p,q}=H_p(H_1(M,{\Z}),H_q^{\rho}(\tilde M)) \Rightarrow H_{p+q}(M,\rho).$$
Recall that 
$M$ has the homotopy type of an $(n+1)$-complex and 
$\tilde M$ has the homotopy type of a bouquet of spheres of dimension $n$, see \ref{univcover}.  

The group $H_q^{\rho}(\tilde M)$ carries the canonical structure
of $H_1(M,{\Z})$-module coming from the corresponding module 
structure on chains. We have the isomorphism:
$H_0^{\rho}(\tilde M)={\C}_{\rho}$  where ${\C}_{\rho}$
is the one dimensional representation of $H_1(M,{\Z})$ 
given by $\rho$. Indeed, if $x$  is a generator
of  $C_0^{\rho}(\tilde M)$ as $H_1(M,{\Z})$-module then
$x-x\cdot g=0$ in $H_0^{\rho}(\tilde M)$. On the 
other hand $x-x\cdot g=x-\rho(g)g^{-1}x$.
Hence $gx=\rho(g)x$ in $H_0^{\rho}(\tilde M)$.

Since for $\rho \ne 1$ one has 
$H_p(H_1(M,{\Z}),{\C}_{\rho})=0$, it follows
from the vanishing theorem in the last section that 
the term $E_2$ has only one 
horizontal row: $q=n$. This yields the claim 1.
We have 
$$H_0(H_1(M,{\Z}),H_n^{\rho}(\tilde M))=
H_n^{\rho}(\tilde M)_{Inv}=
H_n(\tilde M) \otimes_{H_1(M,{\Z})}{\C}_{\rho}$$
On the other hand, taking tensor product with ${\C}_{\rho}$
in the resolution 
$\Phi: \Lambda^s \rightarrow \Lambda^t \rightarrow H_n(\tilde M)
\rightarrow 0$
we obtain the resolution of $H_n(\tilde M) \otimes_{H_1(M,{\ Z})} {\C}_{\rho}$
in which the matrix of $\Phi_{\rho}: \Lambda^s \otimes {\C}_{\rho} \rightarrow 
\Lambda^t \otimes {\C}_{\rho}$ is obtained from the matrix 
of $\Phi$ by replacing its entries by values of the entries at 
$\rho$. Hence if $\rho$ belongs to the set of zeros of the $k$-th 
Fitting ideal (and $k$ is maximal with this property), then the corank of $\Phi_{\rho}$ is $k$. This yields  the second claim.

Let us consider the exact sequence:
\begin{equation}\label{fiveterms}
H_{n+1}(M) \rightarrow H_{n+1}(H_1(M,{\Z})) \rightarrow
H_n(\tilde M) \otimes_{H_1(M,{\Z})} {\C} \rightarrow
H_n(M) \rightarrow H_n(H_1(M,{\Z})) \rightarrow 0
\end{equation}
corresponding to the spectral sequence:
\begin{equation}\label{spectralseq}
H_p(H_1(M,{\Z}),H_q(\tilde M)) \Rightarrow H_{p+q}(M) 
\end{equation}
We have 
\begin{equation}\label{cokernel}
{\rm dim \ Coker} 
(H_{n+1}(M) \rightarrow H_{n+1}(H_1(M,{\Z})))=
{\rm dim \ Ker} H_{n+1}(H_1(M,{\Z}))^* \rightarrow
H_{N+1}(M)^*)
\end{equation}
The latter kernel (using Kronecker pairing 
identification $H_i^*=H^i$ over ${\C}$)
is isomorphic to 
${\rm dim \ Ker} H^{n+1}(H_1(M)) \rightarrow
H^{n+1}(M)$. Since $H_1(M,{\Z})={\Z}^r$ we have 
$\Lambda^i(H_1(M))=H^i(H_1(M))$ with the isomorphism 
provided by the cup product.
Hence the dimension in (\ref{cokernel}) is equal to  
$${\rm dim \ Ker} ( \Lambda^{n+1}(H_1(M)) \rightarrow
H^{n+1}(M)).$$ 
Therefore, using the sequence (\ref{fiveterms})
and the equality ${\rm dim}H_n(H_1(M))={r \choose n}$, we obtain
$${\rm dim} H_n(\tilde M) \otimes_{H_1(M,{\Z})} {\C}=
{\rm dim \ Ker} (\Lambda^{n+1}(H_1(M)) \rightarrow
H^{n+1}(M))+{\rm dim} H_n(M)-{r \choose n}.$$
This yields the last claim in Theorem \ref{locsys}.

\end{proof}

Notice that the claim \ref{locsys}.3 is a generalization of 
a result in \cite{abcov} and that the space of local systems
with non vanishing cohomology in the cases when $Y={\C}^{2}$
and $Y={\C}^{n+1}$ was studied in \cite{Hodge} and \cite{INNC}
respectively.    
The following result is a generalization of Example (6.1.8) in \cite{D},
where $Y$ was assumed to be a smooth germ and the proof uses properties of the vanishing cycle functor and a generalization of Prop. 4.6 from \cite{INNC} 
where the case of  $X$ smooth was treated. 

\begin{cor}\label{monodromy}

Let $F$ be the Milnor fiber of the reduced germ $g:(Y,0) \to (\C,0)$ which defines the divisor $X$ in $Y$. With the above assumptions, the monodromy action 
$h^j:H^j(F,\C) \to H^j(F,\C) $ is trivial for $j \le n-1$.

\end{cor}

\begin{proof} Let $\rho _a$ be the representation sending each elementary loop
to the same complex number $a \in \C^*$. Then it is well-known, see for instance \cite{D}, Corollary 6.4.9, that
$$\dim H^q(M,\rho_a)=\dim \ker(h^q-a Id)+\dim \ker(h^{q-1}-a Id).$$
The unipotence follows by applying this equality to $a \ne 1$.

To obtain triviality of the monodromy action, notice that due to the Milnor's  
fibration theorem, the Milnor fiber $F$ is homotopy equivalent to 
the infinite cyclic cover of $M$. Hence, it is a quotient of the universal 
abelian cover $\tilde M$ by the action of the kernel of 
$\pi: \pi_1(M) \rightarrow {\bf Z}$
where $\pi$ sends each elementary loop to $1$. Let us consider 
the corresponding spectral sequence:
\begin{equation}\label{spectralseq2}
H^p({\rm Ker}\pi,H^q(\tilde M,{\C})) \Rightarrow H^{p+q}(F,{\C})
\end{equation}
for this action of the group ${\rm Ker} \pi={\Z}^{r-1}$ on the universal 
abelian cover $\tilde M$. Notice that this 
is a spectral sequence of ${\C}[t,t^{-1}]$ modules 
where the action on $H^p({\rm Ker}\pi,H^q(\tilde M))$
is the standard action of the generator of $\pi_1/{\rm Ker}\pi$ and the 
action of $t$ on cohomology of the Milnor fiber is the 
monodromy action. 
 Since, by the theorem \ref{univcover}, 
$H^q(\tilde M)=0$ for $1 \le q <n$ we have 
$n-1$ zero-rows in the term $E_2$ and hence the isomorphism
$H^j(F,{\C})=H^j({\rm Ker}\pi,H^0(\tilde M))$ for $1 \le j \le n-1$. 
Since the map of the 
classifying spaces ${S^1}^r \rightarrow S^1$ corresponding to the 
homomorphism $\pi$ has trivial monodromy, the action of 
$\pi_1/{\rm Ker} \pi$ on $H^j({\rm Ker} \pi,{\C})$ is trivial for any $j$ 
in the range $0 \le j \le n-1$ and the claim 
follows. 
\end{proof}

\begin{remark} \rm
One can obtain the triviality of the monodromy action also using 
mixed Hodge theory, at least for $j < n-1$. 
See for details \cite{DS}, Theorem 0.2. 
Note that the above proof shows that 
$\dim H^j(F)= {r-1 \choose j}$ for $j \leq n-1$ (cf. \cite{INNC}). 
\end{remark}

\section{Homology of local systems (higher rank case)}

In this section we work with weaker assumptions on the germs $(Y,0)$ and $(D_j,0)$ as above. Indeed, we simply need that $M$ has only locally complete intersection singularities (which is weaker than asking $(Y,0)$ to be an isolated singularity) and that there is a $\Q$-Cartier divisor, say $D_1$, among the divisors $D_j$ such that the INNC condition for $X$ holds only along $D_1$. In particular, the singularities of the divisors $D_j \setminus D_1$ for $j>1$ can be arbitrary.

To start, note that if $D_1 \setminus \{0\}$ is contained in the smooth part of the space  $Y \setminus \{0\}$, then one has an
elementary loop
$\delta _1$ which goes once about the irreducible divisor $D_1$ (in a transversal at a smooth point). It follows that a rank $m$
local system $\L$ on $Y \setminus X$ which corresponds to a representation
$$\rho : \pi _1(Y \setminus X) \to GL_m(\C)$$
gives rise to a monodromy operators $T_1=\rho( \delta _1)$. Of course, both 
$\delta _1$ and  $T_1$ are well-defined only up-to conjugacy. The following result should be compared to the vanishing part of Theorem 0.2 in \cite{DS}.

\begin{thm} \label{locsys2}

Let $\L$ be a local system on
$M$ such that

(i) $M$ is a locally complete intersection and $D_1$ is an irreducible $\Q$-Cartier divisor, i.e. there is an integer $m$ such that $mD_1$ is the zero set of a holomorphic germ;

(ii) $D_1 \setminus \{0\}$ is contained in the smooth part of the space  $Y \setminus \{0\}$ and $X$ has only normal crossings along $D_1 \setminus \{0\}$;
 
(iii) the corresponding monodromy operator $T_1$ has not 1 as an eigenvalue.

 Then $H^k(M,\L)=0$ for $k<n$ and for $k>n+1.$

\end{thm}

\begin{proof}

For this proof we assume that the (good) representatives for our germs 
$Y, D_j,...$ exist as closed analytic subspaces in an open ball $B$ of radius $2 \epsilon$ centered at the origin. This implies in particular that $Y$ is a Stein space, as well as
$Y \setminus X$, which is the complement of the zero set of a holomorphic function on $Y$. Such a Stein space has the homotopy type of a CW complex of dimension at most $(n+1)$ by \cite{H2}, and this already gives
 $H^k(Y \setminus X,\L)=0$ for $k>n+1.$

These representatives are good in the sense that all the germs $Y, D_j,...$
have a conic structure inside the ball $B$ such that the corresponding retractions are the same for all these germs. We represent the links $L_Y,L_X,L_{D_j},...$ as the intersections of these representatives inside $B$ with a sphere
$S$ of radius $\epsilon$. In such a way we have an inclusion
$i_{\epsilon}:L_Y \to Y$ and a retraction $r_{\epsilon}:Y^* \to L_Y$ , with $Y^*=Y \setminus \{0\}$, which induces inclusions and retractions for the other germs.\\
The main tool for the proof below is the theory of constructible (resp. perverse) sheaves. For all necessary background material on this subject we refer to \cite{KS} and \cite{D}.

Let $i:Y \setminus X \to Y \setminus D_1$ be the inclusion and set $\F^*=Ri_*\L \in
 D^b_c( Y \setminus D_1)$, $\F^*_1=\F^*|(L_Y \setminus L_{ D_1})$.
The constructible sheaf complex $\F^*$ has constant cohomology sheaves along the fibers of the retraction $r_{\epsilon}$ (since the topology is constant along such a fiber). It follows that 
$$H^k(Y \setminus X,\L)=\H^k( Y \setminus D_1,\F^*)=
\H^k( L_Y \setminus L_{D_1},\F^*_1).$$
Let $j_1: L_Y \setminus L_{D_1} \to L_Y$ be the inclusion and note that
$$Rj_{1*}\F^*_1=Rj_{1!}\F^*_1$$
exactly as in \cite {EV} and \cite{CDO}, the key points being the assumptions (ii) and (iii) in the above statement.
Since the link $L_Y$ is compact, it plays the role of the compact algebraic variety in  \cite {EV} and \cite{CDO}, and we get
$$\H^k( L_Y \setminus L_{D_1},\F^*_1)=\H^k_c( L_Y \setminus L_{D_1},\F^*_1)$$
for any integer $k$. \\
The new difficulty we encounter here is that $ L_Y \setminus L_{D_1}$ is not
a Stein space (not even properly homotopically equivalent to a Stein space
as the retraction $r_{\epsilon}: Y \setminus D_1 \to  L_Y \setminus L_{D_1}$
is not proper!), hence the vanishing for the last hypercohomology group is not
obvious.

We proceed as follows. We apply first Poincar\'e-Verdier Duality on the 
real semialgebraic set $ L_Y \setminus L_{D_1}$ and get
$$\H^k_c( L_Y \setminus L_{D_1},\F^*_1)^{\vee}=\H^{-k}( L_Y \setminus L_{D_1},
D_{\R}\F^*_1).$$
Here $D_{\R}\F^*_1$ is the dual sheaf of $\F^*_1$
in this real setting. Note that we can also consider the (complex) dual sheaf
$D\F^* \in  D^b_c( Y \setminus D_1)$. It follows that
$$D\F^*|( L_Y \setminus L_{D_1})=D_{\R}\F^*_1[1]$$
since the inclusion $ L_Y \setminus L_{D_1} \to  Y \setminus D_1$
is normally nonsingular in the sense of \cite{GM2} (this is what corrsponds to a non-characteristic embedding in the sense of \cite{KS} in the case of singular spaces).

This yields the following isomorphisms
$$\H^{-k}( L_Y \setminus L_{D_1},
D_{\R}\F^*_1)=\H^{-k}( L_Y \setminus L_{D_1},D\F^*[-1]|( L_Y \setminus L_{D_1}))= \H^{-k-1}( L_Y \setminus L_{D_1},D\F^*|( L_Y \setminus L_{D_1})). $$
Here we are again in the presence of a constructible sheaf complex, namely
$D\F^*$, whose cohomology sheaves are constant along the fibes of the retraction $r_{\epsilon}$.
This implies that
$$\H^{-k-1}( L_Y \setminus L_{D_1},D\F^*|( L_Y \setminus L_{D_1}))=
\H^{-k-1}( Y \setminus {D_1},D\F^*)=\H^{k+1}_c( Y \setminus {D_1},\F^*)$$
the last isomorphism coming from  Poincar\'e-Verdier Duality on the smooth
algebraic variety $ Y \setminus {D_1}$.

Now it is time to note that the shifted local system $\L[n+1]$ is a perverse sheaf on the locally complete intersection variety
$M$ and hence $\F^*[n+1]$ is a perverse sheaf on the variety $ Y \setminus {D_1}$ since the morphism
$i$ is Stein and quasi-finite. It follows that
$$\H^{k+1}_c( Y \setminus {D_1},\F^*)=\H^{k-n}_c( Y \setminus {D_1},\F^*[n+1])=0$$
 for any $k<n$ by Artin's Vanishing Theorem in the Stein setting, see \cite{KS}.

\end{proof}

\bigskip 

\bibliographystyle{amsalpha}

\begin{flushleft}
{\sc
Alexandru Dimca : 

   Laboratoire Bordelais d'Analyse et G\'eom\'etrie\\
  Universit\'e Bordeaux I\\
33405 Talence Cedex, FRANCE. } 

\end{flushleft}

\begin{flushleft}
{\sc
Anatoly Libgober : 

Department of Mathematics\\
University of Illinois at Chicago,\\
Chicago,IL 60607, USA. } 

\end{flushleft}

\end{document}